\documentclass[14pt, reqno]{amsart}
\usepackage{amssymb}
\usepackage{amsmath,amssymb, amsthm}
\usepackage[colorlinks,urlcolor=red]{hyperref}
\usepackage{graphicx}
\theoremstyle{plain}
\newtheorem{thm}{Theorem}[section]

\newtheorem{lem}[thm]{Lemma}

\numberwithin{equation}{section}\theoremstyle{definition}
 \newtheorem{defn}[thm]{Definition}

 \newtheorem{Question}[thm]{Question}
 \newtheorem{rem}[thm]{Remark}

\allowdisplaybreaks

\numberwithin{equation}{section}
\newcommand{\Z}{\Bbb Z}
\newcommand{\N}{\Bbb N}
\newcommand{\R}{\Bbb R}

\newenvironment{pf}
    {{\noindent  \textrm{\textbf{Proof }}}~~}

\begin{document}

\title[   ]
{ Title :  The Feichtinger Conjecture for Exponentials}
\author[]{ Author :   W. Lawton}
\date{}
\thanks{ $^{*}$Corresponding author: Wayne Lawton  }
\thanks{\it Email address: {\rm matwml@nus.edu.sg }}
\maketitle

\vspace*{-0.5cm}

\begin{center}
{\footnotesize $^{1}$Department of Mathematics,  National University of Singapore
}
\end{center}

\begin{center}
{\footnotesize Address : Block S17, 10 Lower Kent Ridge Road, Singapore 119076
 }
\end{center}

\hrulefill

{\footnotesize \noindent {\bf Abstract.}
The Feichtinger conjecture for exponentials asserts that the following property holds for every fat Cantor subset B of the circle group: the set of restrictions to B of exponential functions can be covered by Riesz sets. In their seminal paper on the Kadison-Singer problem, Bourgain and Tzafriri proved that this property holds if the characteristic function of B has Sobolev regularity. Their probability based proof does not explicitly construct a Riesz cover. They also showed how to construct fat Cantor sets whose characteristic functions have Sobolev regularity. However, these fat Cantor sets are not convenient for numerical calculations. This paper addresses these concerns. It constructs a family of fat Cantor sets, parameterized by their Haar measure, whose characteristic functions have Sobolev regularity and their Fourier transforms are Riesz products. It uses these products to perform computational experiments that suggest that if the measure of one of these fat Cantor sets B is sufficiently close to one, then it may be possible to explicitly construct a Riesz cover for B using the Thue-Morse minimal sequence that arises in symbolic topological dynamics.
\vskip0.5cm



\hrulefill

\section{ Introduction }
\label{S:1}
We let $\N = \{1,2,3,\dots\}, \Z, \R, \mathbb C,$ and $\mathbb T = \R/\Z$  denote the natural numbers, integers, reals, complex numbers and the circle group with Haar measure $\mu.$ Throuought this paper $B$ denotes a Borel subset of $\mathbb T$ with $\mu(B) > 0,$ $F$ denotes a nonempty subset of $\Z,$ and $\chi_B$ and $\chi_F$ denote their characteristic functions. $F$ is an arithmetic set if $F = j+n\Z$ for $j \in \Z, n \in \N.$ We define exponential functions
$e_k(t) = e^{2\pi i k t}, k \in \mathbb Z, t \in \mathbb T;$
$E(F) = \{ \,  e_k \, : \, k \in F \, \};$
$$P(F) = \hbox{trigonometric polynomials spanned by } E(F) \hbox{ whose norm } ||f||_2 = 1;$$
\begin{equation}
\label{eq:1.2}
    \alpha(B,F) =
    \inf \, \{ \, \int_{\, t \in B} \, |f(t)|^2 \, dt \, : \, f \in P(F) \, \}.
\end{equation}
$(B,F)$ is a {\bf Riesz pair} if $\alpha(B,F) > 0,$ or equivalently if $E(F)\chi_B$ is a Riesz basic sequence \cite{CH03}. $\{ F_j  :  j = 1,...,n  \}$
is a {\bf Riesz cover} for $B$ if each $(B,F_j)$ is a Riesz pair and $\cup_{j=1}^{n} F_j = \Z.$ This paper studies the
\\ \\
{\bf Feichtinger Conjecture for Exponentials (FCE)} Every $B$ has a Riesz cover.
\\ \\
The FCE is a special case of the Feichtinger Conjecture (FC): Every unit norm Bessel sequence is a finite union of Riesz basic sequences, which was formulated in (\cite{CA05}, Conjecture 1.1). Casazza and Crandel \cite{CC06} proved that the FC is equivalent to a yes answer to the following problem which has remained open since it was formulated in 1959 \cite{KS59}:
\\ \\
{\bf Kadison-Singer Problem (KSP)} Does every pure state on the C$^*$-algebra $\ell^{\infty}(\mathbb Z)$ admit a unique extension to the C$^*$-algebra of bounded operators $B(\ell^2(\mathbb Z))$?
\\ \\
A {\bf fat Cantor} is a set that is homeomorphic to Cantor's ternary set and has positive Haar measure. Lemma \ref{lem:2.1} shows that the FCE is equivalent to the assertion: every fat Cantor set has a Riesz cover. In their seminal paper on the KSP, Bourgain and Tzafriri proved a result (\cite{BT91}, Theorem 4.1) that implies $B$ has a Riesz cover whenever $\chi_B$ is in the Sobolev space $H^s(\mathbb T)$ for some $s > 0.$ However, their existence proof does not explicitly construct Riesz covers. They also proved a result (\cite{BT91}, Corollary 4.2) that implies $\chi_B \in H^{s}(\mathbb T)$ for all $s < \frac{1}{2}$ whenever $\mathbb T \backslash B = \cup_{j=1}^{\infty} O_j$ where $O_j$ are pairwise disjoint open intervals that satisfy $\mu(O_j) < 2^{-j},$ thus showing the existence of fat Cantor sets that have Riesz covers. This is surprising since Lemma \ref{lem:2.2} shows that if $B$ is a fat Cantor set then $(B,F)$ is not a Riesz pair for a class of sets $F$ that includes the class of arithmetic sets.
\\
This paper has four main results:
\\ \\
(i) Construction of {\bf ternary fat Cantor sets} such that the Fourier transforms of their characteristic functions are Riesz products described by Equation \ref{eqn:3.1}.
\\ \\
(ii) Proof of Theorem \ref{thm:3.1} which shows that ternary fat Cantor sets satisfy $\chi_B \in H^s(\mathbb T)$ for every $s < 1 - \frac{\log 2}{\log 3} \approx 0.3691$ so they have Riesz covers. Ternary fat Cantor sets differ from those constructed by Bourgain and Tzafriri because the lengths of the open intervals $O_j$ removed have algebraic decay $j^{-\log 3/\log 2}$ rather than exponential decay $2^{-j}.$ The proof uses Lemma \ref{lem:2.3} which provides a refinement of the standard Paley-Littlewood decomposition that
Bourgain and Tzafriri used to prove (\cite{BT91}, Corollary 4.2).
\\ \\
(iii) Computation of estimates of $\alpha(B,F),$ where B is a ternary fat Cantor set and $\chi_F = \cdots 10010110.0110100110010110 \cdots$ is the Thue-Morse minimal sequence \cite{TH06}, \cite{MO21}. These estimates suggest that $\{F,1+F,2+F\}$ is a Riesz cover for $B$ if $\mu(B)$ is sufficiently close to $1.$
\\ \\
(iv) Proof of Theorem \ref{thm:3.2} that shows $S(F)$ is convex whenever $\chi_B$ is a minimal sequence.
\\ \\
Results (iii) and (iv) relate the FCE to the field of symbolic dynamics. We
give a brief review of the concepts from this field that we use in this paper.
\\ \\
Let $A$ be any finite set with the discrete topology, equip $A^{\Z}$ with the product topology (it is homeomorphic to Cantor's ternary set). The symbolic dynamical system (over $A$) is the pair  $(A^{\mathbb Z},\sigma),$ where $\sigma$ is the shift homeomorphism defined by
$$(\sigma \, b)(n) = b(n-1), b \in A^{\mathbb Z}, n \in \mathbb Z.$$
A sequence $b \in A^{\Z}$ is {\bf minimal} if its orbit closure
$\overline {\{ \, \sigma^n(b) \, : \, n \in \Z \, \} }$
is a minimal closed shift-invariant set. Zorn's lemma ensures the existence of minimal sequences in any nonempty closed shift invariant subset.
$F \subset \Z$ is {\bf syndetic} if there exists $n \in \N$ such that $\cup_{j=0}^{n-1} (j+F) = \Z,$
{\bf thick} if for every $n \in \N$ there exists $k \in \Z$ such that $k + \{0,2,3,...,n-1\} \subset F,$
and {\bf piecewise syndetic} if $F = F_s \, \cap \, F_t$ where $F_s$ is syndetic and $F_t$ is thick.
Minimal sequences are characterized by a result of Gottschalk \cite{G46}, \cite{GH55} that says a sequence $b$ is minimal if and only if for every finite $G \subset \Z,$ the set
$\{ \, n \in \Z \, : \, \sigma^n(b)|_{G} = b|_{G} \, \}$ is syndetic. Gottschalk's theorem shows that the Thue-Morse sequence is a minimal sequence and it can be also used to construct other minimal sequences. Choosing $G = \{0\}$ implies that if $F$ is nonempty and $\chi_F$ is a minimal sequence in $\{0,1\}^{\Z}$ then $F$ is syndetic so for some integer $n,$ $\cup_{j=0}^{n-1} (j + F) = 0.$ Therefore, if $(B,F)$ is a Riesz pair and $\chi_F$ is a minimal sequence then $B$ has a Riesz cover.
Furstenberg  in (\cite{FU81},Theorem 1.23) used Gottschalk's result for symbolic dynamics over the set $\{1,...,n\}$ to prove that if $\cup_{j=1}^{n} F_j = \mathbb Z$ then one of the sets $F_j$ is piecewise syndetic.
In (\cite{LA10}, Theorem 1.1) we used Furstenberg's result to prove that $B$ has a Riesz cover if and only if there exists nonempty set $F$ such that $(B,F)$ is a Riesz pair and $\chi_F$ is a minimal sequence. This result reduces the construction of a Riesz cover for $B$ to the construction of a single set $F$ such that $(B,F)$ is a Riesz pair and $\chi_F$ is a minimal sequence.
Paulsen \cite{PA08} investigated the relationship between the Kadison-Singer Problem and syndetic sets and in \cite{PA10} he used methods from operator algebras (completely positive maps and multiplicative domains) to independently derive the key results in our paper \cite{LA10} that relate the FCE to syndetic sets.
\section{Preliminaries}
$L^2(\mathbb T)$ and $\ell^2(\mathbb Z)$ are Hilbert spaces with scalar products
$(f,g) = \int_{\, x \in \mathbb T} f(x) \overline {g(x)} d\mu(x)$
and
$(a,b) = \sum_{\, n \in \mathbb Z} a(n) \overline {b(n)}$
and associated norms $|| \ \  ||_2$ and the Fourier transform $L^2(\mathbb T) \ni f \rightarrow \widehat f \in \ell^2(\Z),$ defined by
$$\widehat f(k) = (f,e_k) = \int_{\mathbb T} \, f(x)\, e_k(-x) \, dx,  \ \ k \in Z,$$
is a unitary surjection. Sobolev spaces for $s \geq 0$ are defined by
$$
    H^s(\mathbb T) = \{ \, f \in L^2(\mathbb T) \, : \, \sum_{k \in \Z} |{\widehat f}(k)|^{2} \, |k|^{2s} < \infty \, \}.
$$
A function $f \in L^2(\mathbb T)$ is called Sobolev regular if $f \in H^s(\mathbb T)$ for some $s > 0.$
We observe that $\L^1(\mathbb T) \subset L^2(\mathbb T)$ is a Banach algebra under convolution and each
exponential function $e_k$ defines a {\bf multiplicative linear functional} or {\bf character}
$\ell^1(\mathbb T) \ni f \rightarrow \widehat f(k) \in \mathbb C.$
$C(\mathbb T)$ denotes the Banach space of continuous complex valued functions on $\mathbb T$ with the infinity norm $|| \, ||_{\infty}$ and $M(\mathbb T)$ denotes the Banach algebra of complex valued measures on $\mathbb T$ with the total variation norm and the convolution product. The Riesz Representation Theorem asserts that $M(\mathbb T)$ is the linear dual of $C(\mathbb T).$ For $t \in \mathbb Z,$ $\delta_t$ denotes the Dirac measure at $t$ defined by $\delta_t(f) = f(t), f \in C(\mathbb T).$ A measure $\nu$ is positive if $\nu(f) \geq 0$ whenever $f \geq 0,$ a probability measure if it is positive and $\nu(1) = 1,$ and discrete if it is a countable linear combination of Dirac measures. We identify $L^1(\mathbb T)$ with the set of absolutely continuous measures in $M(\mathbb T).$
The {\bf maximal ideal space} $\widehat M(\mathbb T)$ is the set of {\bf generalized characters} which define multiplicative linear functionals on $M(\mathbb T)$ via the Gelfand correspondence.
For $\nu \in M(\mathbb T)$ and $t \in \mathbb T$ we define translation $\tau_t\nu(f) = \nu(f(\cdot - t))$ and for $n \in \N$ we define the dilation $d_n\nu$ by $d_n\nu(f) = \sum_{k=0}^{n-1} \int_{0}^{1} f(\frac{x+k}{n}) d\nu(x).$ Clearly $d_n\nu$ has period $\frac{1}{n}$ since $\tau_{\frac{1}{n}}d_n \nu = d_n \nu.$
\begin{defn} \label{def:2.1}
The {\bf spectral envelope} of $F$ is the set of probability measures on $\mathbb T$ given by
\begin{equation}\label{eq:3.1}
    S(F) = \mbox{weak}^{*}-\mbox{closure} \, \{ \, |f|^2 \, : \, f \in P(F)  \}.
\end{equation}
\end{defn}
If $j \in \Z, n \in \N$ then $S(j+F) = S(F)$
and $S(nF) = \{ \, d_n\nu \, : \, \nu \in S(F) \, \}.$
If $\nu \in S(F)$ and $t \in \mathbb T$ then $\tau_t \, \nu \in S(F).$
The functions
$f_n = \frac{1}{\sqrt n} (e_1 + \cdots + e_n) \in P(\mathbb Z)$
and the sequence of Fejer kernels $K_n = |f_n|^2$ converges weakly to $\delta_0 \in S(\mathbb Z).$ If $\nu \in M(\mathbb T)$ is a probability measure then $K_n * \nu$ are nonnegative trigonometric polynomials so the Riesz-Fejer spectral factorization theorem implies there exists $Q_n \in P(\mathbb T)$ such that $K_n * \nu = |Q_n|^2.$ Since $K_n*\nu$ converges weakly to $\nu,$ it follows that $\nu \in S(\Z).$ Therefore $S(\Z)$ consists of all probability measures and $S(j + n\Z)$ contains the discrete measure
$d_{\frac{1}{n}}\delta_0 = \sum_{k=0}^{n-1} \delta_{\frac{k}{n}}.$
\\ \\
If $B_0 \subset B$ are Borel sets and $B$ does not have a Riesz cover then $B_0$ does not have a Riesz cover.
Assume that a Borel set $B \subseteq \mathbb T$ contains an interval $[a,b]$ with $b > a.$ Choose $n \in \N$ with $n(b-a) \geq 1$ and set $F_j = j+n\mathbb Z \, : \, j = 0,...,n-1.$ If $f \in P(F_j)$ then $|f|^2$ has period $\frac{1}{n}$ and hence $\alpha(B,F_j) \geq 1/n, \, j = 0,...,n-1.$ Therefore $B$ has a Riesz cover consisting of arithmetic sets. If a Borel set $B_1$ with $\mu(B_1) > 0$ does not have a Riesz cover then, since $\mu$ is inner regular, there exists a closed $B_0 \subset B_1$ with $\mu(B_0) > 0.$ Clearly $B_0$ is nowhere dense.
A theorem of Brouwer (\cite{Kechris}, Theorem 7.4) shows if $\mu(B) > 0$ and $B$ is closed, nowhere dense, and perfect (has no isolated points) then $B$ is a fat Cantor set.
\begin{lem}
\label{lem:2.1}
If a Borel set $B_1$ with $\mu(B_1) > 0$ does not have a Riesz cover then there exists a fat Cantor set $B \subseteq B_1$ such that $B$ does not have a Riesz cover. Therefore, the FCE is equivalent to the assertion that every fat Cantor set has a Riesz cover.
\end{lem}
\begin{pf}
By the preceding argument $B_1$ contains a nowhere dense closed set $B_0 \subset B_1$ with $\mu(B_0) > 0.$ We use transfinite to define a collection of nonincreasing subsets $B_\gamma \subseteq B_0$ indexed by ordinals as follows: $B_{\gamma+} = $ set of limit points of $B_\gamma$ and for limit ordinals $B_\beta = \cap_{\gamma < \beta} B_\gamma.$ The Cantor-Baire Stationary Principle implies that there exists
a countable ordinal $\gamma$ such that $B_{\gamma+} = B_{\gamma}.$ Let $B = B_{\gamma}.$ Then $B$ is perfect. A set of isolated points is countable so has Haar measure zero. Since we remove a countable number of such sets from $B_0$ to obtain $B,$ $\mu(B) > 0.$ Since $B$ is also closed and nowhere dense Brouwer's theorem implies that $B$ is a fat Cantor set. $\Box$
\end{pf}
\begin{lem}\label{lem:2.2}
If $B$ is closed, $\nu \in S(F).$ and $\nu(B) = 0$ then $(B,F)$ is not a Riesz pair. If $B$ is a fat Cantor set and $S(F)$ contains a discrete measure then $(B,F)$ is not a Riesz pair. A fat Cantor set can not have a Riesz cover consisting of arithmetic sets.
\end{lem}
\begin{pf}
If $B$ is closed then the map $S(F) \ni \nu \rightarrow \nu(B)$ is upper semi-continuous since $\limsup_{\nu_j \rightarrow \nu} \nu_{j}(B) \leq \nu(B),$ hence
\begin{equation}
\label{eq:2.3}
    \alpha(B,F) =
    \inf \, \{ \nu(B) \, : \, \nu \in S(F) \, \}
\end{equation}
and this implies the first assertion.
If $\nu \in S(F)$ is discrete then there exists sequences $c_j \geq 0, t_j \in \mathbb T, j \in \N$
such that $\nu = \sum_{j\in \N} c_j \delta_{t_j}.$
The Baire Category Theorem implies there exists $t \notin \cup_{j\in \N} (B - t_j).$ Then $\tau_t \, \nu \in S(F)$ and
$$\tau_t \, \nu(B) = \sum_{j \in \N, t_j+t \in B} c_j = \sum_{j \in \N, t \in (B-t_j)} c_j = 0$$
and this implies the second assertion. The fact that spectral envelopes of arithmetic sets contain discrete measures implies the
third assertion. $\Box$
\end{pf}
\begin{rem}\label{rem:2.1}
In \cite{LA10} we describe Bohr sets and show (\cite{LA10}, Theorem 2.1) that if $F$ is a Bohr set then $S(F)$ contains discrete measures. Therefore, by Lemma \ref{lem:2.2}, if $B$ is a fat Cantor set and $F$ is a Bohr set then $(B,F)$ is not a Riesz pair.
\end{rem}
\begin{lem}\label{lem:2.3}
Let $\sigma_j, j \in \N$ be a sequence of positive real numbers that satisfies the following two conditions:
$$(i) \limsup_{j \rightarrow \infty} \frac{\sigma_{j+1}}{\sigma_{j}} < 1$$
and
$$(ii) \liminf_{j\rightarrow \infty} \frac{\sigma_{j+1}}{\sigma_{j}} > 0.$$
If $c > 0,$ $p > 0$ and $f \in L^2(\mathbb T)$ satisfies
\begin{equation}\label{eqn:2.1}
    ||f(\cdot - \sigma_j) - f(\cdot)||_{2}^{2} \leq c \, \sigma_{j}^{p}, \ \ j \in \N
\end{equation}
then $f \in H^s(\mathbb T)$ for all $s \in (0,\frac{p}{2}).$
\end{lem}
\begin{pf}
Equation \ref{eqn:2.1} implies that $\widehat f$ satisfies
\begin{equation}\label{eqn:2.2}
    \sum_{k \in \Z} |\widehat f(k)|^2 \, \frac{4\sin^2\pi k \sigma_j}{\sigma_{j}^p} \leq c, \ \ j \in \N.
\end{equation}
Condition (ii) ensures that there exists $\theta \in (0,\pi)$ such that
$$\frac{\theta}{\pi}\sigma_{j}  \leq (1-\frac{\theta}{\pi}) \sigma_{j+1}, \  \ j \in \N.$$
We set $c_1 = c(1-\theta/\pi)^{2s}/(4\sin^2 \theta)$ and observe that
\begin{equation}\label{eqn:2.3}
   \sum_{|k| \geq \frac{\theta}{\pi}\sigma_{1}^{-1}} |\widehat f(k)|^2 |k|^{2s} \leq
    \sum_{j \in \N} \ \sum_{\frac{\theta}{\pi} \leq |k|\sigma_j \leq (1-\frac{\theta}{\pi})} \, |\widehat f(k)|^2 \, |k|^{2s} \leq
   c_1 \, \sum_{j \in \N} \, \sigma_{j}^{p-2s}.
\end{equation}
Condition (i) ensures that this sum converges whenever $s \in (0,\frac{p}{2}).$ $\Box$
\end{pf}
%
\section{
Main Results}
\noindent {\bf Construction of Ternary Fat Cantor Sets}
\\ \\
For every $\gamma \in (0,1)$ the following construction gives a fat Cantor set $B \subset \mathbb T$ such that $\mu(B) = \gamma.$ Start with the interval $S_0 = [-\frac{1}{2},\frac{1}{2}]$ and remove the middle open interval having length $\frac{1}{3}(1-\gamma)$ to obtain a set $S_1$ equal to the union of two disjoint equal length closed intervals. From each of these two intervals remove the middle open interval having length $\frac{1}{9}(1-\gamma)$ to obtain a set $S_4$ equal to the union of four disjoint equal length closed intervals. Continue in this manner to construct a decreasing sequence of closed sets $S_j$ each the union of $2^j$ closed intervals having length
$$L_j = \gamma \, 2^{-j} + (1-\gamma) \, 3^{-j}.$$
Construct $S = \cap_{j\in \N} S_j$ and $B = S + Z \subset \mathbb T.$ Clearly $\mu(B) = \gamma$ and $B = -B.$ Let $x_j$ be the distance between the center of the rightmost interval in $S_{j-1}$ and the rightmost interval in $S_j,$ let
$I_j = [-\frac{1}{2}L_j,\frac{1}{2}L_j],$
and define the sequence of discrete measures
$$
\nu_j =
    \left(\frac{1}{2}\delta_{-x_1}+\frac{1}{2}\delta_{x_1}\right)*
\left(\frac{1}{2}\delta_{-x_2}+\frac{1}{2}\delta_{x_2}\right) * \cdots *
\left(\frac{1}{2}\delta_{-x_j}+\frac{1}{2}\delta_{x_j}\right).
$$
Then
$$x_j = \frac{1}{2}(L_{j-1} - L_j), \ \ j \in \N$$
and
$$
    \chi_{S_j} = 2^j \, \chi_{I_j} * \nu_j.
$$
Since we have weakly convergent sequences
$\chi_{S_j} \rightarrow \chi_B$
and
$2^j \, \chi_{I_j} \rightarrow \mu(B) \, \delta_0$
it follows that
$$\mu(B)\, \nu_j \rightarrow \chi_B.$$
Therefore the Fourier transform of $\chi_B$ equals the Riesz product (\cite{ZYG}, Section 7, Chapter 5)
\begin{equation}\label{eqn:3.1}
\widehat {\chi_{B_\beta}}(k) = \mu(B) \prod_{j \in \N} \cos(2\pi x_j k), \ \ k \in \Z.
\end{equation}
Equation \ref{eqn:3.1} provides an efficient method to compute $\widehat \chi_B.$
\begin{thm}\label{thm:3.1}
If $B$ is a ternary fat Cantor set then $\chi_B \in H^s(\mathbb T)$ for $s < 1 - \frac{\log 2}{\log 3}$ so $B$ has a Riesz cover.
\end{thm}
\begin{pf}
Assume that $B$ is a ternary fat Cantor set and set $\gamma = \mu(B).$
Set $\sigma_j = 3^{-j}(1-\gamma), j \in \N$ and $p = 1 - \frac{\log 2}{\log 3}.$
Lemma \ref{lem:2.3} implies that it suffices to show that there exists $c > 0$ such that
\begin{equation}\label{eqn:3.2}
    ||\chi_B(\cdot - \sigma_j) - \chi_B||_{2}^{2} \leq c \, \sigma_{j}^{p}, \ \ j \in \N.
\end{equation}
The Borel subsets of $\mathbb T$ form an abelian group under the Boolean operation
$$B_1 \Delta B_2 = (B_1 \cup B_2) \backslash (B_1 \cap B_2)$$
with identity $\phi,$ $B_1 \Delta B_1 = \phi,$ $\mu(B_1 \Delta B_2) \leq \mu(B_1) + \mu(B_2),$
and $||\chi_{B_1} - \chi_{B_2}||_{2}^{2} = \mu(B_1 \Delta B_2).$
Observe that since $S_{j}$ consists of the union of $2^j$ closed intervals separated by distance $\geq \sigma_j,$
$$\mu((S_{j}+\sigma_j) \Delta S_{j}) \leq 2(1-\gamma)\left(\frac{2}{3}\right)^{j}, \, j \in \N.$$
Furthermore
$$\mu(S_j \Delta B) = \mu(S_j \backslash B) = \sum_{k = j}^{\infty} 2^{k} \sigma_{k+1} =
(1-\gamma) \, \left(\frac{2}{3}\right)^{j}.$$
Inequality \ref{eqn:3.2} holds with $c = 4(1-\gamma)^{1-p}$ since
$(B+\sigma_j) \, \Delta \, B = [(B+\sigma_j) \, \Delta \, (S_j\, + \, \sigma_j)] \, \Delta [B \, \Delta \,  S_j] \, \Delta \, [(S_j+\sigma_j) \, \Delta S_j]$
implies that
$$||\chi_B(\cdot - \sigma_j) - \chi_B||_{2}^{2} \leq 4(1-\gamma) \, \left(\frac{2}{3}\right)^{j}.$$
\end{pf}
\begin{thm}\label{thm:3.2}
If $\chi_F$ is a minimal sequence then $S(E)$ is convex. Furthermore,
\begin{equation}\label{extreme}
    \alpha(B,F) = \inf \{ \, \nu(B) \, : \, \nu \in S_e(F) \, \}
\end{equation}
where $S_e(F)$ is the set of extreme points in $S(E).$
\end{thm}
\begin{pf}
Set $Q(F) = \{ \, |f|^2 \, : \, f \in P(F) \, \}.$ Since $S(E)$ is the weak$^{*}$ closure of $Q(F),$ to prove that $S(F)$ is convex it suffices to show that the convex combination of any two elements in $Q(F)$ is in $S(E).$ Let $f, h \in P(F)$ and let $a, b \in [0,1]$ satisfy $a^2 + b^2 = 1.$
Gottshalk's theorem implies there exists a sequence $n_j \in \N$ converging to $\infty$ with
$e_{n_j}\, h \in P(F), \, j \in \N.$
Define the sequence
$$g_j = \frac{|f + e_{n_j}\, h|^2}{||f + e_{n_j}\, h||_2}, \ \ j \in \N.$$
Then $g_j \in Q(F), \ \ j \in \N$ and the Riemann-Lebesgue lemma implies that
$$
    \lim_{j \rightarrow \infty} g_j = a^2 |f|^2 + b^2 |h|^2
$$
thus proving that $S(F)$ is convex. $S_e(F)$ is nonempty since the the Krein-Milman theorem
implies that $S(F)$ is the weak*-closure of the set of convex combinations of points in $S(F).$
Since $S(F)$ is separable, Choquet's theorem implies that every element $\nu \in S(F)$ is represented
by a probability measure on $S_e(F) \subset S(F),$ from which Equation \ref{extreme} follows. \\
\end{pf}
{\bf Optimization Algorithm to Estimate} $\alpha(B,F).$
We now describe a computational approach to estimate $\alpha(B,F)$ under the assumption that
$\chi_B$ is a minimal sequence. Let $B(\ell^2(\Z))$ denote the $C^{*}$-algebra of bounded operators on the Hilbert space $\ell^2(\Z).$
Define the Laurent operator
$L_B \in B(\ell^2(\Z))$ by the Toeplitz matrix
$[L_B](j,k) = \widehat \chi_B(k-j), \, j,k \in \mathbb Z$
and define the restriction operator
$R_F  \, : \, \ell^2(\Z) \rightarrow \ell^2(F)$ by
$R_F(a)(k) = a(k), \ \ a \in \ell^2(\Z), \, k \in F,$
so the adjoint
$R_{F}^{*} : \ell^2(F) \rightarrow \ell^2(\Z)$
is given by
$$R_{F}^{*}(b)(k) = b(k) \hbox{ if}\,  k \in F, \hbox{ else} = 0.$$
The matrix $[R_F \, L_B \, R_{F}^{*}]$ for the operator
$R_F \, L_B \, R_{F}^{*} \, : \, \ell^2(F) \rightarrow \ell^2(F)$
is a principle submatrix of the matrix $[L_B]$ for the Laurent operator $L_B.$
Then
\begin{equation}
\label{eqn:3.3}
\alpha(B,F) = \inf \hbox{spec} R_F \, L_B \, R_{F}^{*}
\end{equation}
where spec denotes the spectrum of the restricted operator. For finite $G \subset \Z,$
$\alpha(B,G) = \min eig [R_G \, L_B \, R_{G}^{*}]$ since the later matrix is finite. For infinite
$F$ such that $\chi_F$ is a minimal sequence, define $F_n = [0,n] \cap F, \, n \in \N.$ Gottschalk's
theorem implies that for every finite $G \subset F$ there exists $m \in \N$ such that $G \subset F_n$
whenever $n \geq m.$ Since $F_n$ is an increasing sequence of sets the sequence $\alpha(B,F_n)$ is a nonincreasing sequence of nonnegative numbers. This implies the following result which provides an algorithm to approximate $\alpha(B,F).$
\begin{equation}
\label{eqn:3.4}
\alpha(B,F) = \lim_{n \rightarrow \infty} \alpha(B,F_n).
\end{equation}
\\ \\
{\bf Description of Numerical Experiments}
We used Equation \ref{eqn:3.1} to compute $\widehat \chi_B$ for ternary Cantor sets.
Figure 1, Figure 2 shows the values $\widehat \chi_B(k), k = 1 : 4095$ for $\mu(B) = 0.25, 0.75,$ respectively.
We used Equation \ref{eqn:3.4} to estimate $\alpha(B,F),$ where $\chi_F$ is the Thue-Morse minimal sequence, by $\alpha(B,F_{4095}),$ where $F_{4095} = [0,4095] \cap F.$
For $\mu(B) = 0.25,$ the computed value of $\alpha(B,F_{4095})$ is the negative number
$-1.2261 \times 10^{-14}$
due to the fact that the true value is smaller than machine precision.
For $\mu(B) = 0.75$ the computed value of $\alpha(B,F_{4095})$ is $0.085512$ which is
$385$ trillion times machine precision! What explains this difference? We proved in
(\cite{LA10},Corollary 1.1) that if $B(B,F)$ is a Riesz pair then $D^{+}(F) \leq \mu(B)$
where the upper {\bf Beurling density}
\begin{equation}\label{eqn:3.5}
    D^{+}(F) = \lim_{k \rightarrow \infty} \max_{a \in \R} \frac{|F \cap (a,a+k)|}{k}.
\end{equation}
Here $|F \cap (a,a+k)|$ is the cardinality of $F \cap (a,a+k).$ This result was based on a deep
result of Landau (\cite{LAN67}, Theorem 3) in a form discussed by Olevskii and Ulanovskii \cite{OU08}.
Clearly, if $\chi_F$ is the Thue-Morse minimal sequence then $D^{+}(F) = \frac{1}{2},$ so
for $\mu(B) < \frac{1}{2},$ $\alpha(B,F) = 0.$ This means that trigonometric polynomials having frequencies in $F$ can have their squared moduli {\bf localized} on the set $\mathbb T \backslash B.$ The coefficients of the most localized polynomial having frequencies in the finite set $F_n$ are the entries of the eigenvectors corresponding to the eigenvalue $\alpha(B,F_n)$ of the restricted matrix.
It is an open question if this happens for $\mu(B) \geq \frac{1}{2}.$
The function in Figure 1 displays more intermittency than the function in Figure 2 because the gaps are larger. Perhaps this difference in intermittency can be used to explain the immense difference in the $\alpha$ values.
\\ \\
We used Equation \ref{eqn:3.4} to compute $\alpha(B,F_n)$ as a function of $L = \log_{2} n$
for ten ternary Cantor sets $B$ with $\mu(B) \in \{ \, 0.5, 1.5, 2.5, ... , 9.5 \, \}.$
Figure 3 shows the values of $\alpha(B,F_n)$ and Figure 4 shows the values of $\log \, \alpha(B,F_n)$
for each of the ten sets. Both plots show that for $\mu(B) < \frac{1}{2},$
$$\alpha(B,F) = \lim_{L \rightarrow \infty} \alpha(B,F_n) = 0.$$
However, Figure 1 shows that for $\mu(B) > \frac{1}{2},$ $\alpha(B,F_n)$ decreases as a function of $L$ much slower and Figure 2 suggests that for $\mu(B) > \frac{1}{2},$ the sequence may not converge to $0$ because
$\log \alpha(B,F_n)$ appears to be a convex function of $L.$ If this is the case then for $\mu(B) > \frac{1}{2},$ $\alpha(B,F) > 0$ and $\{F, 1+F, 2+F\}$ is a Riesz cover for $B.$
\\ \\
{\bf Suggestions for Further Research} Theorem \ref{thm:3.2} shows that characterization of the set of extreme points $S_e(F)$ in the spectral envelopes of integer subsets $F$ such that $\chi_F$ is a minimal sequence is crucial to understanding the FCE. For such a set $F$ consider the dynamical system $(X(F),\sigma)$ where $X(F)$ is the orbit closure of $\chi_B.$ Then $X(F)$ has at least one shift invariant ergodic probability measure $\zeta.$
\begin{rem}\label{rem:3.1}
If $\chi_F$ is the Thue-Morse minimal sequence then $\zeta$ is unique \cite{KE68}.
\end{rem}
Define $X_1(F) = \{ \, b \in X(F) \, : \, b(0) = 1 \, \}.$ For $g \in L^2(X,\zeta,\sigma)$ define $\sigma g(x) = g(\sigma \, x).$ Then the sequence $(\sigma^j \, g,g)$ is positive definite so by the Herglotz theorem there exists a positive measure $\nu_g \in M(\mathbb T)$ such that
$\nu_g(e_j) = (\sigma^j \, g,g), \ \ j \in \Z.$
Define the set
$$\Sigma(F,\zeta) = \{ \, \nu_g \, : \, g \in L^2(X,\zeta), \hbox{ support}(g) \subseteq X_1(F) \, \}.$$
The Birkhoff ergodic theorem can be used to show that
$
\Sigma(F,\zeta) \subset S(F).
$
This fact, together with the fact that $Q(F)$ contains no extreme points, suggests research to answer the
\begin{Question}\label{qn:3.1}
Is $S_e(F) \subseteq \Sigma(F,\zeta) ?$
\end{Question}
The fact that generalized characters play a crucial role in characterizing the structure of the Banach algebra $M(\mathbb T)$ suggests research to investigate their utility for characterizing spectral envelopes.
%
%
\begin{figure}
\includegraphics[width=0.45\textwidth]{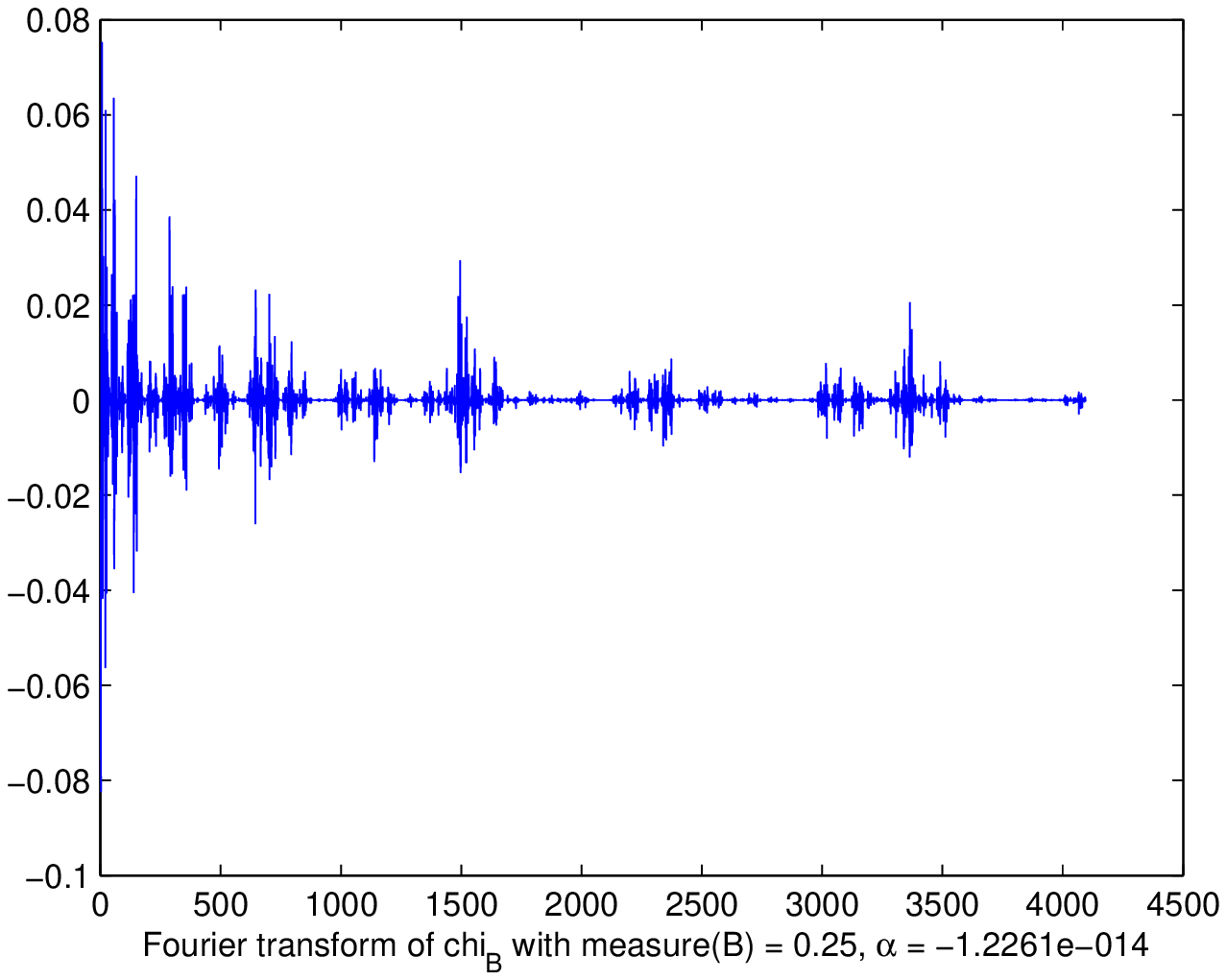}
\includegraphics[width=0.45\textwidth]{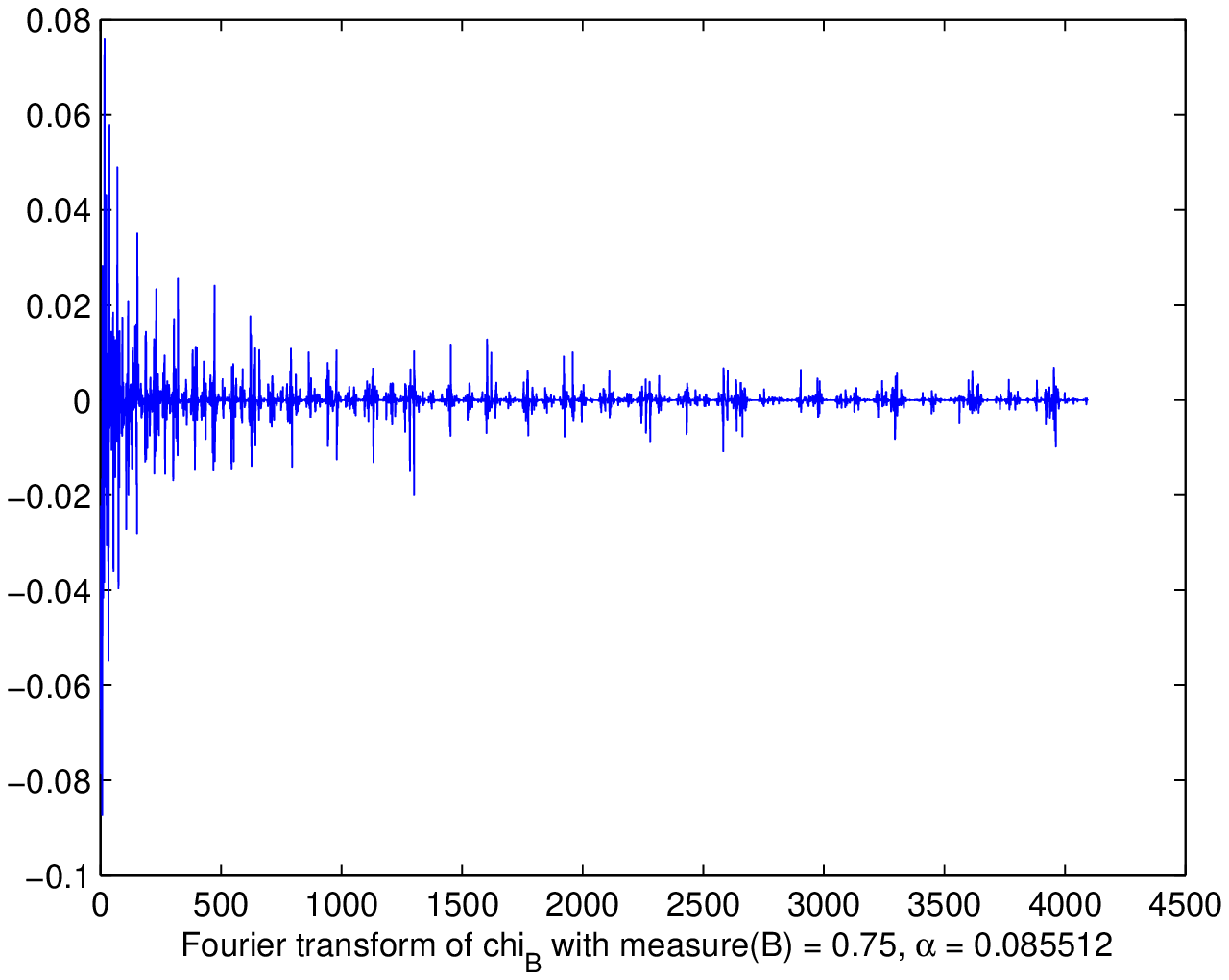}
\caption{Fourier Transform of Characteristic Function of Ternary Fat Cantor Sets}
\end{figure}
\begin{figure}
\includegraphics[width=0.45\textwidth]{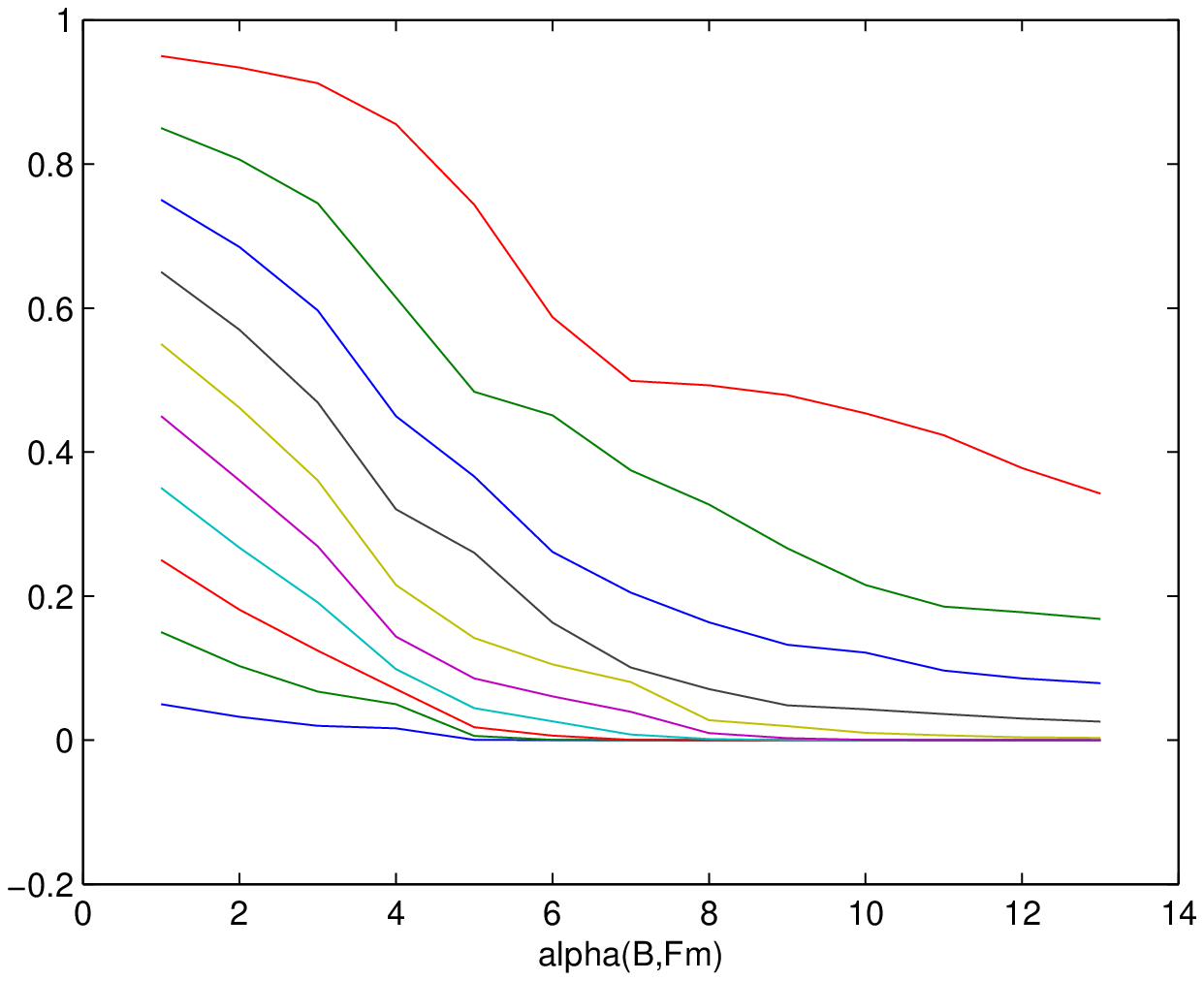}
\includegraphics[width=0.45\textwidth]{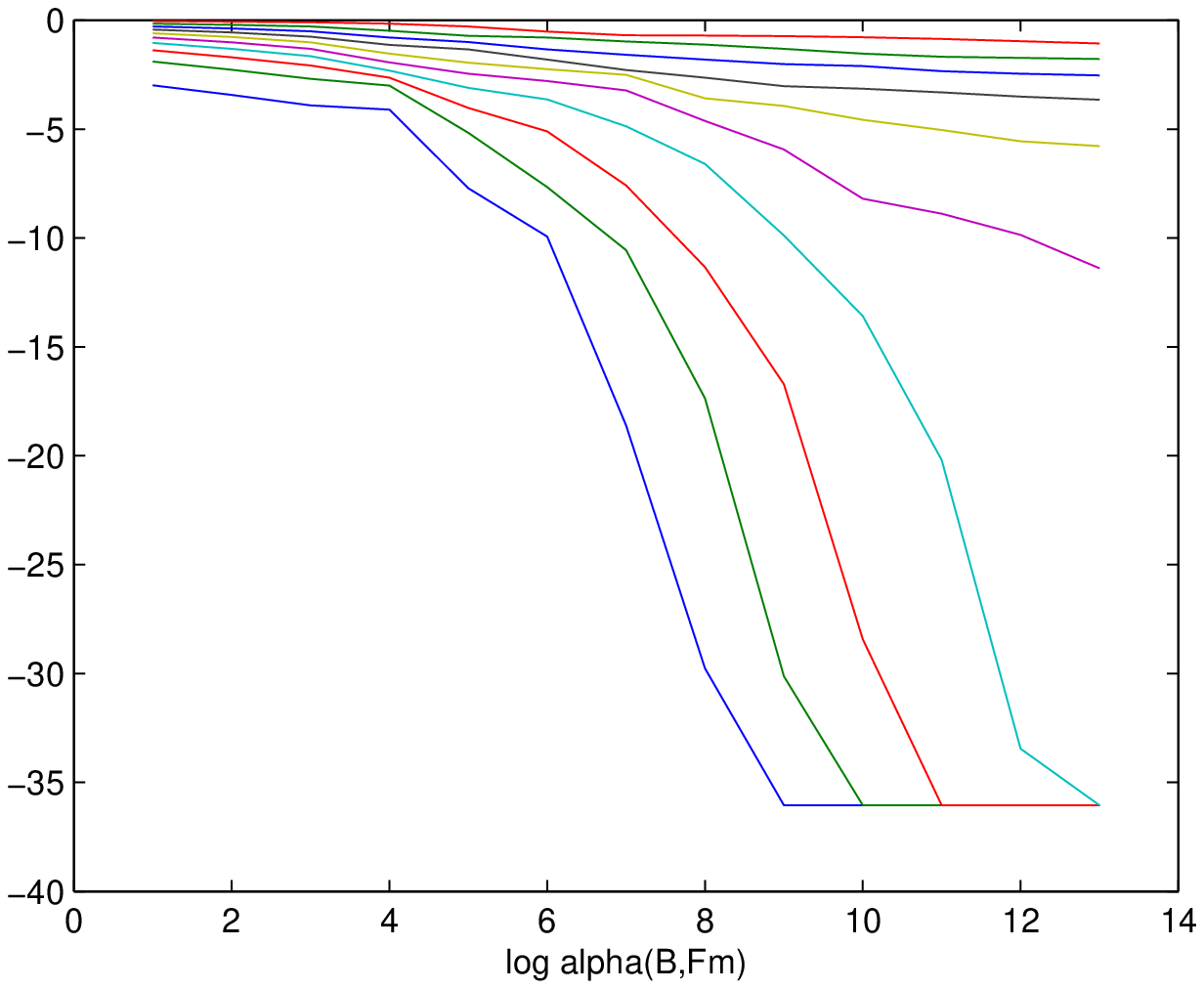}
\caption{alpha(B,F), mu(B) = 0.05 : 0.1 : 0.95, as a function of (1 + log$_2$ size F) }
\end{figure}
\\ \\
\noindent{\bf Acknowledgements} The author would like to thank his colleagues Denny Leung
and Tang Wai Shing for helpful discussions about descriptive set theory and functional analysis.
\bibliographystyle{amsplain}

\end{document}